# Analysis of a Leslie-Gower model with Alle effects, cooperative hunting, and constant placement rates

YONGHUI ZHAO[1]

**ABSTRACT:** This paper investigates the dynamical properties of the Leslie-Gower model with Alle effects, cooperative hunting, and constant placement rates. The conditions for the existence of the triple equilibrium point of the model are first analyzed. Subsequently, the canonical type theory and the qualitative theory of planar systems are applied to obtain that the triple equilibrium point can be a node with a residual dimension of 2 and an equilibrium point with a residual dimension of 3 under different parameter conditions. Finally, it is proved that the system bifurcates with a residual dimension of 2 in the vicinity of the node with cooperative hunting and placement rate as branching parameters.

## 0  Introduction

In ecosystems, there are direct relationships between populations, such as predation, competition, cooperation and parasitism[1]. However, considering the complexity and diversity of populations and their ecological environments, some indirect factors should also be considered between populations, such as the Alle effect, which describes the correlation between population size or density and the fitness of individual populations[2]. In nature, the Alle effect is closely related to the survival and extinction of populations, and in order to protect endangered species, scholars have carried out in-depth research on the Alle effect. In the twentieth century, Stephens[3] and others divided the Allele effect into strong Allele effect and weak Allele effect after research, and Courchamp[4] and others introduced the most typical form of strong Allele effect:

$$f(x) = r\left(1 - \frac{x}{K}\right)(x - m)$$

where $m$ is the threshold for the Allee effect. In order to study the influence of strong Allele effect on predator system, Olivares[5] proposed the following Leslie-Gower model with strong Allele effect:

$$\begin{cases} \dfrac{dx}{dt} = rx\left(1 - \dfrac{x}{K}\right)(x - m) - \lambda xy \\ \dfrac{dy}{dt} = sy\left(1 - \dfrac{y}{cx}\right) \end{cases}, \qquad (0.1)$$

where $x$ and $y$ denote the population densities of the lure and predator, respectively, at time $t$, $r$ is the birth rate of the lure, $K$ is the environmental maximum holding capacity, and $m$ denotes the strong Allee effect. In addition, the stability of the equilibrium point of the (0.1) system, as well as the Hopf branching case, are discussed by Olivares.

In order to effectively capture prey in isolation or in packs, predators cooperate during hunting, e.g., African lions, African wild dogs, birds, ants, etc[6]. Hilker[7] have studied and proposed the following functional response function with cooperative hunting:

$$q(x,y) = (\lambda + ay)x,$$

where $\lambda$ denotes the attack rate of the predator and $a$ denotes the cooperative parameters of the predator in hunting. In order to study the influence of cooperative hunting factors, some scholars added cooperative hunting factors in predator-feeding bait model, for example, penghui [6] et al. considered cooperative hunting factors in Leslie-Gower model with strong Allee effect and built the following model:

$$\begin{cases} \dfrac{dx}{dt} = rx(1-\dfrac{m}{K})(x-m) - (\lambda+ay)xy \\ \dfrac{dy}{dt} = sy(1-\dfrac{y}{cx}) \end{cases}, \quad (0.2)$$

The stability of the equilibrium points of the model as well as the saddle-node and Hopf branching cases are also discussed.

## 1 Model building

In this paper, based on the Leslie-Gower model (0.2) with Alle effect and cooperative hunting, we consider the effect of the placement rate on the system (0.2), and therefore propose the following Leslie-Gower model with Alle effect, cooperative hunting, and a constant placement rate on the prey:

$$\begin{cases} \dfrac{dx}{dt} = rx\left(1-\dfrac{x}{K}\right)(x-m) - (\lambda+ay)xy + h \\ \dfrac{dy}{dt} = sy\left(1-\dfrac{y}{cx}\right) \end{cases}, \quad (1.1)$$

where $x$, $y$ denote the densities of the bait and predator populations, respectively, $r$ the growth rate of the bait, $m$ is the strong Allee effect threshold, $K$ is the maximum environmental holding capacity, $\lambda$ denotes the attack rate of the predator, $a$ denotes the cooperation factor, and $h$ denotes the constant rate of placement on the prey, and denotes the growth rate of the predator. All parameters in model (1.1) are positive, utilizing the following scale transformations:

$$\tilde{x} = \dfrac{x}{K}, \tilde{y} = \dfrac{y}{cK}, \tau = rKt, \tilde{m} = \dfrac{m}{K}, \tilde{\lambda} = \dfrac{c\lambda}{r}, \tilde{a} = \dfrac{c^2 Ka}{r}, \tilde{h} = \dfrac{h}{rK^2}, \tilde{s} = \dfrac{s}{rK}.$$

Replacing $x, y, t, \lambda, a, h, s$ with $\tilde{x}, \tilde{y}, \tau, \tilde{m}, \tilde{\lambda}, \tilde{a}, \tilde{h}, \tilde{s}$ yields the following model:

$$\begin{cases} \dfrac{dx}{dt} = x(1-x)(x-m) - (\lambda+ay)xy + h \\ \dfrac{dy}{dt} = sy\left(1-\dfrac{y}{x}\right) \end{cases}, \quad (1.2)$$

where the threshold for the strong Allee effect takes the range $0 < m < 1$.

## 2 Dynamical states of model equilibrium points

The dynamical state of the equilibrium point can effectively explain and predict the behavior of the ecosystem and provide scientific basis for ecological environmental protection. In the following chapter, the kinetic state of the equilibrium point of the system will be investigated.

In this paper, it is remembered that $A = (m+1-\lambda)^2 - 3(a+1)m$, $B = -(m+1-\lambda)m + 9(a+1)h$,

$a_1 = (\lambda^2 - 2(m+1)\lambda + (m+1)^2 - 3m)/(3m)$, $h_1 = ((m+1-\lambda)m)/(9(a+1))$.

## 2.1 Existence of an equilibrium point

Considering the practical significance of predator-feeder systems, this subsection only studies the existence of positive equilibrium points inside the system (1.2) and obtains the following theorem:

Theorem 2.1: When $0 < m < 1, 0 < \lambda < m+1-\sqrt{3m}, a = a_1, h = h_1$, there exists a triple positive equilibrium of the system $E_1 = (x_1, y_1) = ((m+1-\lambda)/(3(a+1)), (m+1-\lambda)/(3(a+1)))$.

Proof: Solving for the equilibrium point of the system (1.2) means solving the following equation:

$$\begin{cases} x(1-x)(x-m) - (\lambda + ay)xy + h = 0 \\ sy\left(1 - \dfrac{y}{x}\right) = 0 \end{cases}, \tag{2.1}$$

At At that time, $y \neq 0$. From the second equation of equation (2.1), $y = x$, which can be obtained by substituting it into the first equation and organizing it:

$$(1+a)x^3 - (m+1-\lambda)x^2 + mx - h = 0. \tag{2.2}$$

Now let $g(\lambda) = \lambda^2 - 2(m+1)\lambda + (m+1)^2 - 3m$, know by analysis that the function has two positive zeros and is: $\lambda_1 = m+1-\sqrt{3m}$, $\lambda_2 = m+1+\sqrt{3m}$, ,so when $0 < \lambda < m+1-\sqrt{3m}$, $g(\lambda) > 0$. Now consider the case of $0 < m < 1, 0 < \lambda < m+1-\sqrt{3m}, a = g(\lambda)/3m > 0$, $h = ((m+1-\lambda)m)/(9(a+1)) > 0$, to get $A = B = 0$. At this point, combined with Shengjin's formula, we know that the equation (2.2) has a triple root, denoted as $x_1$, where $x_1 = (m+1-\lambda)/(3(a+1))$.

In summary, it can be obtained that at that time $0 < m < 1, 0 < \lambda < m+1-\sqrt{3m}, a = a_1, h = h_1$, the system (1.2) has a triple positive equilibrium point $E_1 = (x_1, y_1)$, ,where $y_1 = x_1$.. Proof.

Note: Since the dynamical properties of the positive equilibrium point under other parameter conditions as well as the bifurcation are similar to those discussed in the literature[8][10], they are not discussed in this paper.

## 2.2 Stability of the equilibrium point

For determining the local stability situation at the equilibrium point of system (1.2), it is generally necessary to go to the eigenvalue situation of the Jacobi matrix corresponding to the equilibrium point of system (1.2). For system (1.2) the Jacobi matrix at the equilibrium point $E_1 = (x_1, y_1)$ is:

$$J|_{E_1} = \begin{bmatrix} -(a+3)x_1^2 + (2m+2-\lambda)x_1 - m & -2ax_1^2 - \lambda x_1 \\ s & -s \end{bmatrix}, \tag{2.3}$$

Based on the eigenvalue analysis of the Jacobi matrix, the following theorem can be obtained:

**Lemma 2.1:** When $0 < m < (11-\sqrt{105})/4$, there is $m+1-(\sqrt{30m})/2 > 0$ holds. Furthermore, when $0 < m < (11-\sqrt{105})/4, 0 < \lambda < m+1-(\sqrt{30m}/2)$, there is $a_1 > 3/2$ holds.

**Theorem 2.2:** Let $s_1 = x_1(2ax_1 + \lambda)$, for the internal triple equilibrium point $E_1 = (x_1, y_1)$ of the system there are the following cases:

(1) When $0 < m < 1, 0 < \lambda < m+1-\sqrt{3m}, a = a_1, h = h_1, s \neq s_1$, $E_1$ is a node of residual dimension 2.

(2) When $0 < m < (11-\sqrt{105})/4, 0 < \lambda < m+1-(\sqrt{30m}/2), a = a_1, h = h_1, s = s_1$, $E_1$ is an equilibrium point of residual dimension 3.

PROOF: Now $J(E_1)$ sets the trace and determinant to $tr(J(E_1))$ and $det(J(E_1))$, respectively, according to (2.3): $Tr(J|_{E(x_1, y_1)}) = -(a+3)x_1^2 + (2m+2-\lambda)x_1 - m - s = 2ax_1^2 + \lambda x_1 - s = s_1 - s$ as well as $Det(J|_{E(x_1, y_1)}) = s(3(1+a)x_1^2 - 2(m+1-\lambda)x_1 + m) = 0$. The following transformation is made to the system (1.2):

$$\begin{cases} u = x - x_1 \\ v = y - y_1 \end{cases},$$

Translating $E_1$ to the origin and Taylor expanding it at the origin yields the following system:

$$\begin{cases} \dfrac{du}{dt} = a_{10}u + a_{01}v + a_{20}u^2 + a_{11}uv + a_{02}v^2 + a_{30}u^3 + a_{12}uv^2 \\ \dfrac{dv}{dt} = b_{10}u + b_{01}v + b_{20}u^2 + b_{11}uv + b_{02}v^2 + b_{30}u^3 + b_{21}u^2v + b_{12}uv^2 + b_{40}u^4 + b_{31}u^3v \\ \quad + b_{22}u^2v^2 + O(|u,v|^5) \end{cases} \quad (2.4)$$

Among them: $a_{10} = x_1(2ax_1 + \lambda) = s_1, a_{01} = -x_1(2ax_1 + \lambda) = -s_1, a_{20} = m+1-3x_1,$

$a_{11} = -\lambda - 2ax_1, a_{30} = -1, a_{12} = -a, b_{10} = s, b_{01} = -s, b_{20} = -\dfrac{s}{x_1}, b_{11} = \dfrac{2s}{x_1}, b_{02} = -\dfrac{s}{x_1},$

$b_{30} = \dfrac{s}{x_1^2}, b_{21} = -\dfrac{2s}{x_1^2}, b_{12} = \dfrac{s}{x_1^2}, b_{40} = -\dfrac{s}{x_1^3}, b_{31} = \dfrac{2s}{x_1^3}, b_{22} = -\dfrac{s}{x_1^3}.$

(a) When $s \neq s_1$, there is $Tr(J|_{E(x_1, y_1)}) \neq 0$ as well as $Det(J|_{E(x_1, y_1)}) = 0$. At this point the matrix $J|_{E(x_1, y_1)}$ J has a zero eigenvalue and non-zero eigenvalues.

In order to change the linearized matrix of the system (2.4) at the origin to the standard form, the following transformation is done:

$$\begin{pmatrix} u \\ v \end{pmatrix} = \begin{pmatrix} 1 & s_1 \\ 1 & s \end{pmatrix} \begin{pmatrix} u_1 \\ v_1 \end{pmatrix},$$

Get:

$$\begin{cases} \dfrac{du_1}{dt} = c_{20}u_1^2 + c_{11}u_1v_1 + c_{02}v_1^2 + c_{30}u_1^3 + c_{21}u_1^2v_1 + c_{12}u_1v_1^2 + c_{03}v_1^3 + O(|u_1,v_1|^4) \\ \dfrac{dv_1}{dt} = (s_1-s)y + d_{20}u_1^2 + d_{11}u_1v_1 + d_{02}v_1^2 + d_{30}u_1^3 + d_{21}u_1^2v_1 + d_{12}u_1v_1^2 + d_{03}v_1^3 \\ \qquad\quad + O(|u_1,v_1|^4) \end{cases}, \quad (2.5)$$

Among them: $c_{20}=0$, $c_{11}=\dfrac{-s(2s_1a_{20}+a_{11}(s+s_1)+2sa_{02})+s_1(2s_2b_{20}+b_{11}(s+s_1)+2sb_{02})}{s_1-s}$,

$c_{02}=\dfrac{-s(s_1^2a_{20}+s_1sa_{11}+s^2a_{02})+s_1(s_1^2b_{20}+s_1sb_{11}+s^2b_{02})}{s_1-s}$, $c_{30}=\dfrac{s(a+1)}{s_1-s}$,

$c_{21}=\dfrac{-s(3s_1a_{30}+2sa_{12}+s_1a_{12})+s_1(3s_1b_{30}+sb_{21}+2s_1b_{21}+s_1b_{12}+2sa_{12})}{s_1-s}$,

$c_{12}=\dfrac{-s(3s_1^2a_{30}+2ss_1a_{12}+s^2a_{12})+s_1(3s_1^2b_{30}+2ss_1b_{21}+s^2b_{21}+2ss_1b_{12}+s^2b_{12})}{s_1-s}$,

$c_{03}=\dfrac{-s(s_1^3a_{30}+s^2s_1a_{12})+s_1(s_1^3b_{30}+s_1^2sb_{21}+s^2s_1b_{12})}{s_1-s}$, $d_{20}=0$,

$d_{11}=\dfrac{2s_1a_{20}+a_{11}(s+s_1)+2sa_{02}-(2s_1b_{20}+b_{11}(s+s_1)+2sb_{02})}{s_1-s}$,

$d_{02}=\dfrac{s_1^2a_{20}+s_1sa_{11}+s^2a_{02}-(s_1^2b_{20}+s_1sb_{11}+s^2b_{02})}{s_1-s}$, $d_{30}=\dfrac{(a_{30}+a_{12})+s_1(b_{30}+b_{21}+b_{12})}{s_1-s}$,

$d_{21}=\dfrac{s_1a_{30}+2sa_{12}+s_1a_{12}-(s_1b_{30}+sb_{21}+2s_1b_{21}+s_1b_{12}+2sa_{12})}{s_1-s}$,

$d_{12}=\dfrac{3s_1^2a_{30}+2ss_1a_{12}+s^2a_{12}-(3s_1^2b_{30}+2ss_1b_{21}+s^2b_{21}+2ss_1b_{12}+s^2b_{12})}{s_1-s}$,

$d_{03}=\dfrac{s_1^3a_{30}+s^2s_1a_{12}-s_1(s_1^3b_{30}+s_1^2sb_{21}+s^2s_1b_{12})}{s_1-s}$.

Now consider the case $s<s_1$ and do the time transformation $dt=(s_1-s)d\tau$ for the system (2.5) to get:

$$\begin{cases} \dfrac{du_1}{d\tau} = e_{11}u_1v_1 + e_{02}v_1^2 + e_{30}u_1^3 + e_{21}u_1^2v_1 + e_{12}u_1v_1^2 + e_{03}v_1^3 + O(|u_1,v_1|^4) \\ \dfrac{dv_1}{d\tau} = v_1 + f_{11}u_1v_1 + f_{02}v_1^2 + f_{30}u_1^3 + f_{21}u_1^2v_1 + f_{12}u_1v_1^2 + f_{03}v_1^3 + O(|u_1,v_1|^4) \end{cases}. \quad (2.6)$$

Let the center of the system (2.6) be prevalent as: $v_1 = h(u_1) = \sigma_1 u_1^2 + \sigma_2 u_1^3 + O(|u_1|^4)$, then there is:

$$(2\sigma_1 u_1 + 3\sigma_2 u_1^2 + \cdots)(e_{11}u_1h(u_1) + e_{02}h^2(u_1) + \cdots) - (\sigma_1 u_1^2 + \sigma_2 u_1^3 + O(|u_1|^4)) - f_{11}u_1(\sigma_1 u_1^2 + \sigma_2 u_1^3 + +O(|u_1|^4)) - f_{02}h^2(u_1) - f_{30}u_1^3 - f_{21}u_1^2h(u_1) + \cdots = 0,$$

This is obtained by calculation: $\sigma_1 = 0, \sigma_2 = -f_{30}$, Then the center prevalence is

$v_1 = h(u_1) = -f_{30}u_1^3 + O(|u_1|^4)$, It follows that the The system (2.6) is equivalent to the following system near the origin:

$$\frac{du_1}{d\tau} = e_{30}u_1^3 - e_{11}f_{30}u_1^4 + O(|u_1|^5), \qquad (2.7)$$

included among these:

$$e_{30} = (s(a+1))/(s_1 - s)^2 > 0, \quad e_{11}f_{30} = (s^2 s_1(a+1)(s_1 - s)(4ax_1 + \lambda))/(s_1 - s)^4 \neq 0.$$

due to $e_{30} > 0$, Combined with the literature[11], it is known that $E_1$ is an unstable node of residual dimension 2 when $s < s_1$, and the same analysis shows that $E_1$ is a stable node of residual dimension 2 when $s > s_1$.

(b) When $s = s_1$, there is $Tr(J|_{E(x_1,y_1)}) = 0$ as well as $Det(J|_{E(x_1,y_1)}) = 0$, when the matrix $J|_{E(x_1,y_1)}$ has two zero eigenvalues. In order to change the linearized matrix of the system (2.4) at the origin into standard form, the following transformation is done:

$$\begin{pmatrix} u \\ v \end{pmatrix} = \begin{pmatrix} -s_1 & 0 \\ -s_1 & 1 \end{pmatrix} \begin{pmatrix} U_1 \\ V_1 \end{pmatrix},$$

The following system is obtained:

$$\begin{cases} \dfrac{dU_1}{dt} = V_1 + g_{11}U_1V_1 + g_{02}V_1^2 + g_{30}U_1^3 + g_{21}U_1^2V_1 + g_{12}U_1V_1^2 + O(|U_1,V_1|^4) \\ \dfrac{dV_1}{dt} = h_{11}U_1V_1 + h_{02}V_1^2 + h_{30}U_1^3 + h_{21}U_1^2V_1 + h_{12}U_1V_1^2 + O(|U_1,V_1|^4) \end{cases}, \qquad (2.8)$$

Among them: $g_{11} = a_{11} + 2a_{02}, g_{02} = -a_{02}/s_1, g_{30} = s_1^2(a_{30} + a_{12}), g_{21} = -2s_1 a_{12}, g_{12} = a_{12}$,

$h_{11} = s_1(a_{11} + 2a_{02} - b_{11} - 2b_{02}), h_{02} = b_{02} - a_{02}, h_{30} = s_1^3(a_{30} + a_{12} - b_{30} - b_{21} - b_{12}), h_{21} = s_1^2(-2a_{12} + b_{21} + 2b_{12})$,
$h_{12} = s_2(a_{12} - b_{12})$.

In order to simplify the second order terms in system (2.8), the following transformations are made to system (2.8):

$$\begin{pmatrix} U_1 \\ V_1 \end{pmatrix} = \begin{pmatrix} U_2 \\ V_2 \end{pmatrix} + \begin{pmatrix} \dfrac{g_{11} + h_{02}}{2}U_2^2 + g_{02}U_2V_2 \\ g_{02}U_2V_2 \end{pmatrix}$$

The following system is obtained:

$$\begin{cases} \dfrac{dU_2}{dt} = V_2 + i_{30}U_2^3 + i_{21}U_2^2V_2 + i_{12}U_2V_2^2 + i_{03}V_2^3 + O(|U_2,V_2|^4) \\ \dfrac{dV_2}{dt} = j_{11}U_2V_2 + j_{30}U_2^3 + j_{21}U_2^2V_2 + j_{12}U_2V_2^2 + j_{03}V_2^3 + O(|U_2,V_2|^4) \end{cases}, \qquad (2.9)$$

Among them: $i_{30} = h_{30}, i_{21} = g_{21} + \dfrac{g_{11}(g_{11} + h_{02})}{2} - h_{02}(g_{11} + h_{02}) + \cdots, i_{12} = g_{12} + g_{11}g_{02} + \cdots$,

$i_{03} = g_{02}(1 - h_{02}), j_{11} = h_{11}, j_{30} = h_{30}, j_{21} = h_{21} + \dfrac{h_{11}(g_{11} + h_{02})}{2} + g_{11}h_{02} - h_{02}h_{11}$,

$j_{12} = h_{12} + g_{11}h_{02} + 2g_{02}h_{02} + \cdots, j_{03} = h_{02} - h_{02}^2$.

By computing $j_{11}j_{30} = s_1^4(\lambda + 4ax_1)(1 + (s/x_1^2)) \neq 0$, the system (2.9) is equivalent by the

lemma in the literature [12]:

$$\begin{cases} \dfrac{dU_2}{dt} = V_2 \\ \dfrac{dV_2}{dt} = j_{11}U_2V_2 + j_{30}U_2^3 + (j_{21}+3i_{30})U_2^2V_2 + O(|U_2,V_2|^4) \end{cases}, \quad (2.10)$$

Among them: $j_{11} = -s_1(m+1+4a_1x_1) < 0$, $j_{30} = -s_1^3(a_1+1) < 0$,

$j_{21} + 3i_{30} = (s_1^3 + s_1^2\lambda x_1 + 2s_1a_1\lambda x_1^3 + 2a_1^2 s_1 x_1^4 + 2s_1^2 x_1^2(2a_1-3))/(2x_1^2)$.

Combining this with Lemma 2.1, when $0 < m < (11-\sqrt{105})/4, 0 < \lambda < m+1-(\sqrt{30m}/2)$, $a = a_1$, $h = h_1$,, there are $j_{11} \neq 0$, $j_{30} \neq 0$ and, $j_{21} + 3i_{30} \neq 0$ then an equilibrium point with residual dimension 3. The proof is complete.

## 2.3 Bifurcation studies

By Theorem 2.2, $E_1$ is a node of coset dimension 2 when $0 < m < 1, 0 < \lambda < m+1-\sqrt{3m}$, $a = a_1$, $h = h_1, s \neq s_1$. Consider whether a bifurcation of the system with coset dimension 2 occurs near the equilibrium point $E_1$ when the parameter $(a,h)$ is transformed in a small domain around $(a_1,h_1)$.

**Theorem 3.1:** For the system (1.2), the parameters $a$ and $h$ branching parameters are chosen such that the system undergoes a bifurcation with a residual dimension of 2 near the equilibrium point E1 when the $(a,h)$ parameters are transformed in a small field near $(a_1,h_1)$.

Proof: by Theorem 2.2 when $0 < m < 1, 0 < \lambda < m+1-\sqrt{3m}, a = a_1, h = h_1, s \neq s_1$, the resulting system (2.4) is equivalent to the following system near the origin when the system (1.2) is translated to the origin at the point of flat $E_1$ equilibrium:

$$\frac{du_1}{d\tau} = e_{30}u_1^3 - e_{11}f_{30}u_1^4 + O(|u_1|^5), \quad (3.1)$$

Where $e_{30} > 0$, now make $d\tau = e_{30}dt$ to get the following system:

$$\frac{du_1}{dt} = u_1^3 - \frac{e_{11}f_{30}}{e_{30}}u_1^4 + O(|u_1|^5), \quad (3.2)$$

From the literature [11], the system (3.2) has the following universal open fold:

$$\frac{du_1}{dt} = \eta_1 + \eta_2 u_1 + u_1^3 + O(|u_1|^4) \quad (3.3)$$

This tells us that the system (1.2) bifurcates at the equilibrium point $E_1$ with a residual dimension of 2. The proof is complete.

Note: In addition, a B-T bifurcation with a residual dimension of 3 may occur in the system (1.2) near the equilibrium point, which will not be discussed in this paper due to the limited length of the article, but will be discussed subsequently in other articles.

## 2.4 Conclusion

The article focuses on the dynamical properties of the Leslie-Gower model with the Alle effect, cooperative hunting, and a constant placement rate on prey. First, the existence case of the triple equilibrium point E1 within the model is analyzed. Secondly, it is proved that $E_1$ can be a node with residual dimension 2 and an equilibrium point with residual dimension 3 under different circumstances using regularity theory and theorems such as central prevalence. Finally, using the bifurcation theory, it is analyzed to obtain that for the case of $s \neq s_1$, the system uses the placement rate $h$ of the arresters and the cooperation parameter $a$ as the branching parameter, and the bifurcation of the node with a residual dimension of 2 occurs in the vicinity of the equilibrium point. Meanwhile, the theoretical analysis reveals that after the introduction of the factor of placement rate, the system may have a node with a residual dimension of 2 and an equilibrium point with a residual dimension of 3, while this type of equilibrium point is not observed in the original model. This also indicates that the inclusion of the placement rate on the prey has a significant effect on the dynamical state of the model.

[1]College of Mathematics and Statistics, Chongqing University, Chongqing, 401331, P. R. China

Email address: 1529320906@qq.com